\documentclass[10pt]{amsart}
\usepackage{amssymb,amsbsy,amsmath,amsfonts,amssymb}
\usepackage{latexsym,euscript,exscale}

\title{Finitely approximable groups and actions\\Part I: The Ribes--Zalesski\u\i{} property}
\author {Christian Rosendal}
\address{Department of Mathematics, Statistics, and Computer Science (M/C 249)\\
University of Illinois at Chicago\\
851 S. Morgan St.\\
Chicago, IL 60607-7045\\
USA}
\email{rosendal.math@gmail.com}
\urladdr{http://www.math.uic.edu/$~$rosendal}

\newcommand {\Z}{\mathbb Z}

\newcommand{\normal}{\trianglelefteq}

\newcommand{\tom} {\emptyset}

\newcommand{\til}{\rightarrow}

\newcommand {\del}{ \; \big| \;}

\newcommand{\inv}{^{-1}}

\newtheorem{thm}{Theorem}

\newtheorem{lemme}[thm]{Lemma}

\newtheorem{prop} [thm] {Proposition}
\newtheorem{defi} [thm] {Definition}

\usepackage{fouriernc}

\begin{document}

\thanks{The author was partially supported by NSF grants DMS 0901405 and DMS 0919700. The author is also grateful  for the helpful suggestions of the anonymous referee.}


\maketitle
\begin{abstract}
We investigate extensions of S. Solecki's theorem on closing off finite partial isometries of metric spaces \cite{solecki1} and obtain the following exact equivalence: any action of a discrete group $\Gamma$ by isometries of a metric space is finitely approximable if and only if any product of finitely generated subgroups of $\Gamma$ is closed in the profinite topology on $\Gamma$.
\end{abstract}

\section{Introduction}
Suppose a discrete group $G$ acts by isometries on a metric space $(X,d_X)$. The main question we shall consider here is the conditions on the group $G$ that ensure that such an action can always be finitely approximated. So, before we go any further, we need to precisely state our notion of approximation.

Let $G$ be a group acting on sets $X$ and $Y$ and let $A\subseteq X$ and $F\subseteq G$ be arbitrary subsets. An $F$-{\em map} from $A$ to $Y$ is a function $\pi \colon A\til Y$ such that whenever $g\in F$ and  $x, gx\in A$, then $\pi (gx)=g\pi (x)$. Moreover, if $X$ and $Y$ are relational structures of the same type, an $F$-{\em embedding} is simply an injective $F$-map that is an isomorphism of $A$ with its image $F[A]$.

\begin{defi}
Let $G$ be a group acting by isometries on a metric space $(X,d_X)$. We say that the action is {\em finitely approximable} if for any finite $A\subseteq X$ and $F\subseteq G$ there is a finite metric space $(Y,d_Y)$, on which $G$ acts by isometries, and an isometric $F$-embedding $\pi\colon A\til Y$.
\end{defi}

We note that this is a very strong notion of finite approximation in the sense that we require $\pi$ to be an isometric embedding and not just a map with small distortion. However, as we shall see later, this is not a real restriction.

The goal here is to isolate the properties of a group $G$ that ensure that every isometric action of $G$ is finitely approximable and, in particular, show this is so for finitely generated Abelian groups. These results form the basis of our investigation in our companion paper \cite{generics}. The first result along these lines is due to S. Solecki \cite{solecki1}, who proved that this is verified when $G$ is a finitely generated free group. Solecki's proof was based on earlier work of B. Herwig and D. Lascar \cite{herwig}  on equations in free groups and ultimately relied on the solution to Rhodes' Type II Conjecture proved independently by C. J. Ash \cite{ash} and L. Ribes and P. A. Zalesski\u\i  \;\cite{ribes}. The proof by Ribes and Zalesski\u\i{} used a connection between the Type II Conjecture and the profinite topology on free groups established by J.-E. Pin and C. Reutenauer \cite{pin}, and, indeed, the main result of Ribes and Zalesski\u\i{} is that products of finitely generated subgroups of free groups are closed in the profinite topology. While part of the result presented here is not entirely novel and can be proved with the methods of the literature, it is not clear from Solecki's and Herwig and Lascar's papers that this topological property   of $G$ is exactly what is needed for the approximation result.  Our main goal here is to establish the following exact correspondence between the two properties.
\begin{thm}
The following two properties are equivalent for a countable discrete group $G$.
\begin{enumerate}
\item If $H_1,\ldots, H_n$ are finitely generated subgroups of $G$, their product $H_1H_2\cdots H_n$ is a closed subset in the profinite topology on $G$,
\item any action of $G$ by isometries on a metric space is finitely approximable.
\end{enumerate} 
\end{thm}

\subsection{The profinite topology and finitely approximable groups}
We recall that if $G$ is a discrete group, the {\em profinite topology} on $G$ is the group topology on $G$ generated by the basic open sets
$$
gK,
$$
where $g\in G$ and $K$ is a finite index subgroup of $G$. Since any finite index subgroup $K\leqslant G$ contains a further subgroup $N\leqslant K$, which is both of finite index and normal in $G$, in the definition of the basis, one can always assume that $K$ is moreover normal in $G$. Thus, a subset $S\subseteq G$ is closed in the profinite topology on $G$ if for any $g\in G\setminus S$, there is a finite index (normal) subgroup $K\leqslant G$ such that $gK\cap S=\emptyset$ or, equivalently, such  that $g\notin SK$. Since this is a group topology, i.e., the group operations are continuous, $G$ is Hausdorff if and only if $\{1\}$ is closed, i.e., if for any $g\neq 1$ there is a finite index subgroup $K$ not containing $g$. In other words, $G$ is Hausdorff if and only if it is residually finite.

A stronger notion than residual finiteness is subgroup separability or being LERF (locally extended residually finite). Here a group $G$ is {\em subgroup separable}, or {\em LERF}, if any finitely generated subgroup $H\leqslant G$ is closed in the profinite topology on $G$. Since $\{1\}$, is finitely generated, subgroup separability implies residual finiteness. M. Hall \cite{hall1,hall2} originally proved that free groups are subgroup separable.

However, the even stronger notion that concerns us here is what we shall call the  Ribes--Zalesski\u\i{} property, or property (RZ) for brevity. Namely, a group $G$ is said to have the {\em  Ribes--Zalesski\u\i{} property} if any product $H_1H_2\cdots H_n$ of finitely generated subgroups $H_i\leqslant G$ is closed in the profinite topology on $G$. This property was originally proven for free groups in \cite{ribes} and T. Coulbois \cite{coulbois} showed that if both $G$ and $F$ have property (RZ), then so does $G*F$.


\subsection{Finitely approximable actions} The requirement in the definition of finitely approximable actions that the map $\pi$ be an isometry may seem very strong, but can actually be weakened somewhat without effect. 
\begin{defi} 
Suppose $(X,d_X)$ and $(Y,d_Y)$ are metric spaces and $K\geqslant 1$. A function $f\colon X\til Y$ is said to be a {\em map with constant $K$}  if for some constant $c>0$ and all $x,y\in X$,
$$
cd_X(x,y)\leqslant d_Y(f(x),f(y))\leqslant Kcd_X(x,y).
$$
The minimum $C$ of the numbers $K$ for which $f$ is a map with constant $K$ is said to be the {\em distortion} of $f$. 
In this case, we say that $(X,d_X)$ embeds into $(Y,d_Y)$ with {\em distortion $C$}.
\end{defi}
Note that if we rescale the metric $d_Y$ by a factor $\frac 1{Kc}$, the function $f$ will be $1$-Lipschitz and $f\inv$ will be $K$-Lipschitz on its domain $f[X]$.

A finite set $D$ of non-negative real numbers is said to be a {\em good value set} if $0\in D$ and for all $s,t\in D$, if $s+t\leqslant \max D$, then also $s+t\in D$.
With this definition in hand, the following lemma is trivial to verify.

\begin{lemme}\label{good value}
Suppose $(X,d_X)$ is a metric space and $D$ is a good value set. Then, 
$$
\partial_X(x,y)=
\begin{cases}
\min (s\in D\del d_X(x,y)\leqslant s)        &\text{if $d_X(x,y)\leqslant \max D$,}\\
\max D                                                     &\text{otherwise.}
\end{cases}
$$
defines a metric on $X$. Moreover, if $G$ is a group acting by isometries of $(X,d_X)$, then the same action is also an isometric action on $(X,\partial_X)$. 
\end{lemme}

\begin{defi} If $(X,d_X)$ is a metric space and $A\subseteq X$ is a non-empty finite subset, the {\em expanded distance set} of $A$ is defined by 
$$
{\rm Ex}(A)=\{r_1+\ldots+r_n\del r_i\in d_X[A\times A]\;\&\; r_1+\ldots+r_n\leqslant  {\rm diam}(A)\}.
$$
\end{defi}
Note that ${\rm Ex}(A)$ is a good value set in the above sense.  So with these preliminaries off our hands, we can now show the following equivalence.
\begin{prop}\label{distortion}
Suppose $G$ is a group acting by isometries on a metric space $(X,d_X)$ and that $A\subseteq X$ and $F\subseteq G$ are finite subsets. Then the following are equivalent,
\begin{enumerate}
\item there is a finite ${\rm Ex}(A)$-valued metric space space $(Y,d_Y)$, on which $G$ acts by isometries, and an isometric $F$-embedding $\pi\colon A\til Y$,
\item there is a finite metric space space $(Y,d_Y)$, on which $G$ acts by isometries, and an isometric $F$-embedding $\pi\colon A\til Y$,
\item for every $K>1$, there is a finite metric space space $(Y,d_Y)$, on which $G$ acts by isometries, and an $F$-map $\pi\colon A\til Y$ of distortion at most $K$.
\end{enumerate}
\end{prop}

\begin{proof}
The only non-trivial part is that (3) implies (1). So pick $K>1$ small enough such that for any $s<t$ in ${\rm Ex}(A)$, $Ks<t$. Assume $(Y,d_Y)$ is a finite metric space on which $G$ acts by isometries and that $\pi\colon A\til Y$ is an $F$-map with constant $K$. Also, by rescaling the metric $d_Y$, we can suppose that $\pi$ is $1$-Lipschitz and $\pi\inv$ is $K$-Lipschitz on its domain. By Lemma \ref{good value}, we can define a $G$-invariant, ${\rm Ex}(A)$-valued metric  $\partial_Y$ on $Y$ by 
$$
\partial_Y(x,y)=
\begin{cases}
\min (s\in {\rm Ex}(A)\del d_Y(x,y)\leqslant s),       &\text{if $d_Y(x,y)\leqslant \max {\rm Ex}(A)$;}\\
\max {\rm Ex}(A),                                                     &\text{otherwise.}
\end{cases}
$$
Then, for all $x,y\in A$,
$$
\frac 1Kd_X(x,y)\leqslant d_Y(\pi(x),\pi(y))\leqslant d_X(x,y),
$$
and so, by assumption on $K$, we have for all $x,y\in A$, $\partial_Y(\pi(x),\pi(y))=d_X(x,y)$. It follows that $\pi$ is an isometric $F$-embedding of $(A,d_X)$ into $(Y,\partial_Y)$.
\end{proof}
Thus, approximating with arbitrarily small distortion or isometrically are equivalent and we shall therefore stick to the stronger notion throughout the rest of the paper.

\section{Groups with the Ribes--Zalesski\u\i{} property}

We are now ready for our main equivalence.

\begin{thm}\label{metric}
The following properties are equivalent for a countable group $G$.
\begin{itemize}
\item[(a)] $G$ has property (RZ),
\item[(b)] any action of $G$ by isometries on a metric space is finitely approximable.
\end{itemize}
\end{thm}
We note also that it follows from Theorem \ref{metric} together with  Proposition \ref{distortion}, that if $G$ has property (RZ) and acts by isometries on a rational valued metric space, then any finite approximation can be done by rational valued metric spaces. This will be important for our later applications.

\begin{proof}
Suppose first $G$ is a group with property (RZ) acting  by isometries on a metric space $X$ and $A\subseteq X$ and $F\subseteq G$ are finite subsets. Replacing $d_X$ with the metric $\partial_X$ defined from the good value set ${\rm Ex}(A)$, as in Lemma \ref{good value}, we can, without changing distances between points of $A$, assume that  $d_X$ only takes values in the finite set ${\rm Ex}(A)$.

Now let $Ga_1,\ldots,Ga_m$ list the orbits of $G$ on $X$ that intersect $A$. Without loss of generality, we can suppose that $X=Ga_1\cup\ldots\cup Ga_m$. Let also $N_i=\{g\in G\del ga_i=a_i\}$ and let $G$ act by left-translation on the space of left cosets
$$
G/N_1\sqcup \ldots\sqcup G/N_m.
$$
Then the map
\begin{align*}
\pi \colon X&\til G/N_1\sqcup \ldots\sqcup G/N_m\\
ga_i\in Ga_i&\mapsto gN_i\in G/N_i
\end{align*}
is a conjugacy of $G$-actions and we may therefore assume that $X=G/N_1\sqcup \ldots\sqcup G/N_m$, i.e., that $G$ acts by left-translation on $G/N_1\sqcup \ldots\sqcup G/N_m$ preserving a metric $d$ taking values in the finite set ${\rm Ex}(A)$.

Now, since $A\subseteq G/N_1\sqcup \ldots\sqcup G/N_m$ is finite,  there is some finite set $C\subseteq G$ such that $A\subseteq \{gN_i\del i\leqslant  m\;\&\; g\in C\}$. By enlarging $F$, we can also assume that $1\in F$.

Define finitely generated subgroups $M_i\leqslant  N_i$ by
$$
M_i=\langle g\inv fh\del f\in F,\; g,h\in C\;\&\; g\inv fh\in N_i\rangle
$$
and subsets $E^{ij}_r\subseteq G$ for $i,j\leqslant  m$ and $r\in {\rm Ex}(A)\setminus \{0\}$ by
$$
E^{ij}_r=M_i\{g\inv h\del g,h\in C\; \&\; d(gN_i,hN_j)=r\}M_j.
$$
Note that $E^{ij}_r=(E^{ji}_r)\inv$ and clearly also
$E^{ij}_r=M_iE^{ij}_rM_j$.

We claim that whenever $i_0,\ldots, i_n\leqslant  m$, $r_j,r\in  {\rm Ex}(A)\setminus \{0\}$ and $r_1+\ldots+r_n<r$, we have
$$
1\notin E_{r_1}^{i_0i_1}E_{r_2}^{i_1i_2}\cdots E_{r_n}^{i_{n-1}i_n}E^{i_ni_0}_r.
$$
To see this, suppose toward a contradiction that $h_l\in E^{i_{l-1}i_l}_{r_l}$ and $h\in E^{i_ni_0}_r$ are such that $1=h_1h_2\cdots h_nh$. Then
$$
d(N_{i_0},h_1N_{i_1})=r_1,\; d(N_{i_1},h_2N_{i_2})=r_2,\;\ldots\;, d(N_{i_{n-1}},h_nN_{i_n})=r_n
$$
and 
$$
d(N_{i_n},hN_{i_0})=r
$$
whereby
\begin{align*}
r&=d(N_{i_0},h\inv N_{i_n})\\
&=d(N_{i_0},h_1h_2\cdots h_nN_{i_n})\\
&\leqslant  d(N_{i_0},h_1N_{i_1})+d(h_1N_{i_1},h_1h_2N_{i_2})+\ldots+d(h_1\cdots h_{n-1}N_{i_{n-1}},h_1\cdots h_{n-1}h_nN_{i_n})\\
&=d(N_{i_0},h_1N_{i_1})+d(N_{i_1},h_2N_{i_2})+\ldots+d(N_{i_{n-1}},h_nN_{i_n})\\
&=r_1+r_2+\ldots+r_n<r,
\end{align*}
which is a contradiction. So the claim holds.

We remark that each of the sets $E^{ij}_r$ can be written as a finite union of sets, $M_ig\inv hM_j$ where $g,h\in C$, and these sets are themselves products of left-cosets of finitely generated subgroups of $G$.
Thus, any set
$$
E_{r_1}^{i_0i_1}E_{r_2}^{i_1i_2}\cdots E_{r_n}^{i_{n-1}i_n}E^{i_ni_0}_r
$$
is again a finite union of products of left-cosets of finitely generated subgroups of $G$. Using that for any subgroups $H_i\leqslant G$ and any $g_i\in G$
\begin{align*}
&g_1H_1g_2H_2\cdots g_kH_k\\
&=g_1\cdots g_k\cdot(g_2\cdots g_k)\inv H_1(g_2\cdots g_k)\cdot(g_3\cdots g_k)\inv H_2(g_3\cdots g_k)\cdots g_k\inv H_{k-1}g_k\cdot H_k,
\end{align*}
we see that such sets are actually finite unions of left-translates of products of finitely generated subgroups of $G$. Since $G$ has property (RZ), it follows that the sets 
$$
E_{r_1}^{i_0i_1}E_{r_2}^{i_1i_2}\cdots E_{r_n}^{i_{n-1}i_n}E^{i_ni_0}_r
$$
are closed in the profinite topology on $G$.

Using that the sets are closed and do not contain $1$, we can find a finite index, normal subgroup $K\normal G$ such that whenever $i_0,\ldots, i_n\leqslant  m$ and $r_j,r\in  {\rm Ex}(A)\setminus \{0\}$ satisfy $r_1+\ldots+r_n<r$, we have
$$
1K\cap  E_{r_1}^{i_0i_1}E_{r_2}^{i_1i_2}\cdots E_{r_n}^{i_{n-1}i_n}E^{i_ni_0}_r=\tom,
$$
i.e., as $K$ is normal, 
$$
1\notin  E_{r_1}^{i_0i_1}E_{r_2}^{i_1i_2}\cdots E_{r_n}^{i_{n-1}i_n}E^{i_ni_0}_rK= \big(E_{r_1}^{i_0i_1}K\big)\big(E_{r_2}^{i_1i_2}K\big)\cdots \big(E_{r_n}^{i_{n-1}i_n}K\big)\big(E^{i_ni_0}_rK\big).
$$
We set $H_i=M_iK$ and $B^{ij}_r=E^{ij}_rK$.  Then, as $K$ is normal in $G$,
$$
B^{ij}_r=E^{ij}_rK=KE^{ij}_r=K\inv (E^{ji}_r)\inv=(E^{ji}_rK)\inv=(B^{ji}_r)\inv
$$
and
$$
H_iB^{ij}_rH_j=(M_iK)(E^{ij}_rK)(M_jK)=M_iE^{ij}_rM_jK=E^{ij}_rK=B^{ij}_r.
$$
Finally, whenever $i_0,\ldots, i_n\leqslant  m$ and $r_j,r\in  {\rm Ex}(A)\setminus \{0\}$ satisfy $r_1+\ldots+r_n<r$, we have
$$
1\notin B_{r_1}^{i_0i_1}B_{r_2}^{i_1i_2}\cdots B_{r_n}^{i_{n-1}i_n}B^{i_ni_0}_r.
$$

Using the sets $B^{ij}_r$, we can define an ${\rm Ex}(A)$-valued metric $\rho$ on the left-coset space
$$
G/H_1\sqcup G/H_2\sqcup\ldots\sqcup G/H_m
$$
by setting  for distinct $gH_i, fH_j$,
\begin{displaymath}\begin{split}
\rho(gH_i,fH_j)=\min\big(\max  {\rm Ex}(A),&\\
 \inf\{r_1+\ldots+r_n\del& g\inv f\in B_{r_1}^{ii_1}B_{r_2}^{i_1i_2}\cdots B_{r_n}^{i_{n-1}j} \text{ for some } i_1,\ldots,i_{n-1}\leqslant  m\}\big).
\end{split}\end{displaymath}
Note first that since $H_iB_r^{il}=B_r^{il}$ and $B^{lj}_rH_j=B^{lj}_r$ for any $r\in  {\rm Ex}(A)\setminus \{0\}$ and $l\leqslant m$, this definition does not depend on the choice of representatives $g$ and $f$ from the cosets $gH_i$ and $fH_j$. Also, $\rho$ is easily seen to satisfy the triangle inequality and the positivity condition, so $\rho$ is indeed a metric on 
$$
G/H_1\sqcup G/H_2\sqcup\ldots\sqcup G/H_m.
$$
Moreover, for all $g,f,h \in G$ and $i,j$, we have
$$
\rho(gH_i,fH_j)=\rho(hgH_i,hfH_j),
$$
so $\rho$ is invariant under the action of $G$ by left-translation on $G/H_1\sqcup \ldots\sqcup G/H_m$.

We claim that if $g\inv h\in B^{ij}_r$, then $\rho(gH_i,hH_j)=r$. For if not, we would have $ g\inv h\in B_{r_1}^{ii_1}B_{r_2}^{i_1i_2}\cdots B_{r_n}^{i_{n-1}j}$  for some  $i_1,\ldots,i_{n-1}\leqslant  m$ and $r_1+\ldots+r_n<r$. But, as $h\inv g\in B_r^{ji}$, this implies
$$
1=g\inv hh\inv g\in B_{r_1}^{ii_1}B_{r_2}^{i_1i_2}\cdots B_{r_n}^{i_{n-1}j}B^{ji}_r,
$$
contradicting the assumptions on the sets $B^{ij}_r$.

Now define 
$$
\pi \colon A\til G/H_1\sqcup G/H_2\sqcup\ldots\sqcup G/H_m
$$ 
by $\pi (gN_i)=gH_i$ for all $g\in C$ such that $gN_i\in A$. We claim that $\pi $ is an isometric $F$-embedding of $A$ into $G/H_1\sqcup\ldots\sqcup G/H_m$.

First, to see that $\pi$ is well-defined, i.e., that $\pi(gN_i)$ does not depend on the representative $g$ of $gN_i$, suppose that $g,h\in C$ and $gN_i=hN_i\in A$. Then $g\inv h\in N_i$, so, as $1\in F$, also $g\inv 1h\in M_i\leqslant M_iK= H_i$, whereby
$$
\pi(gN_i)=gH_i=hH_i=\pi(hN_i).
$$

Secondly, to see that $\pi $ is an isometric embedding of $A$ into $Y$, suppose $gN_i,hN_j\in A$ are distinct, where $g,h\in C$. Let $r=d(gN_i,hN_j)>0$, whereby $g\inv h\in E^{ij}_r\subseteq B^{ij}_r$. It therefore follows that $\rho(gH_i,hH_j)=r$. So $\pi$ is an isometry.

Finally, to see that $\pi $ is a partial conjugation of the action of $G$, suppose $gN_i\in A$ and $f\cdot gN_i=hN_i\in A$ for some $g,h\in C$ and $f\in F$. Then $h\inv fg\in N_i$ and so $h\inv fg\in M_i\leqslant  H_i$. Therefore, $fgH_i=hH_i$ and so 
$$
f\pi (gN_i)=fgH_i=hH_i=\pi (hN_i)=\pi(fgN_i).
$$
This shows that if $G$ has property (RZ), then any action of $G$ by isometries is finitely approximable.

\

Suppose now conversely that any isometric action of $G$ is finitely approximable and let $H_1,\ldots, H_n$ be finitely generated subgroups of $G$. We must show that for any $g\notin H_1\cdots H_n$ there is a finite index subgroup $K\leqslant G$ such that $g\notin H_1\cdots H_nK$. So let $H_1,\ldots, H_n$ and $g$ be given. 
We set
$$
X=G/H_1\sqcup G/H_2\sqcup\ldots G/H_n
$$
and define a graph $\Gamma=(X,E)$ with vertex set $X$ and edges $\{fH_i,kH_{i+1}\}$, where $fH_i\cap kH_{i+1}\neq \tom$.
Now set
$$
d_X(fH_i,kH_j)=\min \big(n,\text{ number of edges in the shortest path from $fH_i$ to $kH_j$ in $\Gamma$ }\big).
$$
So $(X,d_X)$ is an integer valued metric space. Also, the left shift action of $G$ on $X$ clearly preserves the distance $d_X$. 

Notice that for any $f,k\in G$, 
$$
d_X(fH_1,kH_n)=\begin{cases}      n-1   & \text{if there is a path $fH_1,f_2H_2,\ldots,f_{n-1}H_{n-1},kH_n$ in $\Gamma$,}\\
                                             n                                 & \text{otherwise.} \end{cases}
$$
Note that if $f_1H_1,f_2H_2,\ldots,f_nH_n$ is a path in $\Gamma$, then 
$$
f_1\inv f_2\in H_1H_2\;\;\;\&\;\;\; f_2\inv f_3\in H_2H_3\;\;\;\&\;\;\ldots\;\;\&\;\;\; f_{n-1}\inv f_n\in H_{n-1}H_n,
$$
whence 
$$
f_1\inv f_n=f_1\inv f_2\cdot f_2\inv f_3 \cdots f_{n-1}\inv f_n\in H_1H_2\cdot H_2H_3\cdots H_{n-1}H_n=H_1H_2\cdots H_n.
$$
Since $g=1\inv g\notin H_1H_2\cdots H_n$, it follows that $d_X(H_1,gH_n)=n$.

Let now $A=\{H_1,H_2,\ldots,H_n, gH_n\}$ and set
$$
F=\{\text{generators of }H_1,H_2,\ldots,H_n\}\cup \{1,g\}.
$$
By the assumption on $G$, there is a finite metric space $(Y,d_Y)$ upon which $G$ acts by isometries and an isometric $F$-embedding $\pi\colon A\til Y$ of $A$ into $Y$. Let now $K_i$ be the stabiliser of $\pi(H_i)$ in $G$. Since $Y$ is finite, $K_i$ is a finite index subgroup of $G$. We claim that $H_i\leqslant K_i$. For if $h$ is a generator of $H_i$, then,  as $h\in F$ and  $H_i,hH_i\in A$, we must have $\pi(H_i)=\pi (hH_i)=h\pi(H_i)$ and so $h\in K_i$. Thus, $K_i$ contains the generators of $H_i$ and hence $H_i\leqslant K_i$.

Now suppose $k_i\in K_i$. Since $H_1,\ldots, H_n\in A$ and $d_X(H_i,H_{i+1})=1$, we have 
\begin{align*}
d_Y\big(k_1\cdots k_i\pi(H_i),k_1\cdots k_ik_{i+1}\pi(H_{i+1})\big)=&d_Y\big(k_1\cdots k_i\pi(H_i),k_1\cdots k_i\pi(H_{i+1})\big)\\
=&d_Y\big(\pi(H_i),\pi(H_{i+1})\big)\\
=&1.
\end{align*}
So, by the triangle inequality,
$$
d_Y\big(\pi(H_1),k_1k_2\cdots k_n\pi(H_n)\big)=d_Y\big(k_1\pi(H_1),k_1k_2\cdots k_n\pi(H_n)\big)\leqslant n-1.
$$
In other words, for any $f\in K_1\cdots K_n$, we have $d_Y(\pi(H_1),f\pi(H_n))\leqslant n-1$. 

Now, $g\in F$ and $H_n, gH_n\in A$, so $g\pi(H_n)=\pi(gH_n)$, whence 
$$
d_Y\big(\pi(H_1),g\pi(H_n)\big)=d_Y\big(\pi(H_1),\pi(gH_n)\big)=d_X(H_1,gH_n)=n.
$$
So it follows that  $g\notin K_1K_2\cdots K_n$, whence $g\notin H_1H_2\cdots H_nK_n\subseteq K_1\cdots K_n$.
\end{proof}

We should make a few comments on the above proof. Namely, one  sees that in order for a product $H_1\cdots H_n$ of finitely generated subgroups of $G$ to be closed in the profinite topology, it suffices that any action of $G$ by isometries on a metric space with distance set $\{0,1,\ldots,n\}$ is finitely approximable. Conversely, with somewhat more care, one can check that if any 
product $H_1\cdots H_n$ of finitely generated subgroups of $G$ is closed in the profinite topology, then any action of $G$ by isometries on a metric space with distance set $\{0,1,\ldots,n\}$ is finitely approximable. 

Another type of actions of interest to us is actions by automorphisms on a graph. Here a {\em graph} is a pair $(X,E_X)$ of a non-empty set $X$ and a symmetric, irreflexive relation $E_X$ on $X$. So an action $G\curvearrowright X$ of a discrete group by automorphisms of a graph $(X, E_X)$ is said to be {\em finitely approximable} if for any finite $A\subseteq X$ and $F\subseteq G$ there is a finite graph $(Y, E_Y)$ on which $G$ acts by automorphisms and an $F$-embedding $\pi\colon A\til Y$ embedding $(A, E_X)$ into $(Y,E_Y)$.

Similarly, an action $G\curvearrowright X$ of a discrete group by permutations on a set $X$ is {\em finitely approximable} if for any finite $A\subseteq X$ and $F\subseteq G$ there is a finite set $Y$ on which $G$ acts by permutations and an injective $F$-map $\pi\colon A\til Y$. Then, identifying a discrete set with a metric space with distance set $\{0,1\}$ and a graph with a metric space with distance set $\{0,1,2\}$, the proof above show the following probably well-known facts.
\begin{prop}Let $G$ be a discrete group.
\begin{enumerate}
\item [(A)] The following are equivalent.
\begin{itemize}
\item Any action of $G$ by automorphisms of a graph is finitely approximable,
\item if $H$ and $K$ are finitely generated subgroups of $G$, then $HK$ is closed in the profinite topology on $G$.
\end{itemize}
\item [(B)] The following are equivalent.
\begin{itemize}
\item Any action of $G$ by permutations of a set is finitely approximable,
\item $G$ is subgroup separable.
\end{itemize}
\end{enumerate}
\end{prop}

The {\em Rado} or {\em random graph} is the countable graph $\bf R$, i.e., a set with a symmetric irreflexive relation, defined up to isomorphism by the following extension property: if $X$ is a finite graph and $\phi\colon Y\til \bf R$ is an embedding of an induced subgraph, then $\phi$ entends to an embedding of $X$ into $\bf R$. 

E. Hrushovski investigated the case of graphs in \cite{hrushovski}, and showed that any finite graph can be extended to a larger finite graph in such a way that any partial automorphism of the former extends to a full automorphism of the latter. Modulo the existence of an ultrahomogeneous universal graph, namely the random graph, this is equivalent to showing that any action of a  free group by automorphisms of a graph is finitely approximable. So, by the above proposition, this implies the special case of the Ribes--Zalesski\u\i{} Theorem that a product of two finitely generated subgroups of a free group is closed in the profinite topology.


\section{Locally finite dense subgroups of automorphism groups}
The {\em rational Urysohn metric space} ${\mathbb Q\mathbb U}$ is the countable metric space all of whose distances are rational and satisfying the following extension property:
if  $X$ is a finite metric space with rational distances and $\phi\colon Y\til {\mathbb Q\mathbb U}$ is an isometric embedding of a subset $Y\subseteq X$, then $\phi$ extends to an isometric embedding of $X$ into ${\mathbb Q\mathbb U}$.

Using the extension property, it is easy to see that any countable metric space with rational distances embeds isometrically into ${\mathbb Q\mathbb U}$ and also, by a standard back and forth argument, the above properties define ${\mathbb Q\mathbb U}$ uniquely up to isometry.

The rational Urysohn space is useful for us in that it and its isometry group function as a universal framework for certain arguments about metric spaces. To see this, we need to introduce a few notions. 
\begin{defi}
An isometry $f\colon A\til B$ between finite subsets $A$ and $B$ of the rational Urysohn space ${\mathbb Q\mathbb U}$ is said to be a {\em finite partial isometry} of ${\mathbb Q\mathbb U}$.
\end{defi}
So the restriction of any full isometry of ${\mathbb Q\mathbb U}$ to a finite subset is a finite partial isometry, but more importantly, by a back and forth argument, any finite partial isometry of ${\mathbb Q\mathbb U}$ extends to a full isometry of ${\mathbb Q\mathbb U}$. 

In fact, we also have the following interesting fact due to V. V. Uspenski\u\i{} \cite{uspenskii}; namely, if $G$ is a group acting by isometries on a finite subspace $A\subseteq {\mathbb Q\mathbb U}$, then the action of $G$ extends to an action by isometries on all of ${\mathbb Q\mathbb U}$. To see this, we can without loss of generality suppose that $G$ is finite. Also, modulo an inductive construction, it suffices to show that for any one-point extension $B\supseteq A$, there is a further finite extension $C\supseteq B$ and an action of $G$ on $C$ extending the action of $G$ on $A$.  We identify the unique point in $B\setminus A$ with $1\in G$ and extend the metric $d$ on $B=A\sqcup\{1\}$ to all of $C=A\sqcup G$, by letting $d(x,h)=d(h\inv x, 1)$ for $x\in A$ and $h\in G$ and setting
$$
d(g,h)=\min \big(d(g,x)+d(x,h)\del x\in A\big)
$$
for all $g,h\in G$.
This is easily seen to be a rational-valued metric extending the metric on $B$ and the invariance under the left-shift action by $G$ is trivial.

We equip the group ${\rm Isom}({\mathbb Q\mathbb U})$ of isometries of ${\mathbb Q\mathbb U}$ with the {\em permutation group topology}, that is, the basic open neighbourhoods of $f\in {\rm Isom}({\mathbb Q\mathbb U})$ are of the form 
$$
U(f,x_1,\ldots,x_n)=\{g\in {\rm Isom}({{\mathbb Q\mathbb U}})\del g(x_i)=f(x_i),\; i\leqslant n\},
$$
where $x_1,\ldots, x_n\in {\mathbb Q\mathbb U}$. Since ${\mathbb Q\mathbb U}$ is countable, it is easy to see that ${\rm Isom}({\mathbb Q\mathbb U})$ is a {\em Polish} group, i.e., a separable and completely metrisable topological group.

We shall now use the correspondence developed above in conjunction with the main result of Coulbois \cite{coulbois} to give an alternative proof of an unpublished result of Solecki (we include the result here with Solecki's permission).

\begin{thm}[S. Solecki]
The isometry group ${\rm Isom}({\mathbb Q\mathbb U})$ of the rational Urysohn metric space has a dense, locally finite subgroup.
\end{thm}

\begin{proof}
Let $f_0,f_1,f_2,\ldots$ be a listing of a dense subset of ${\rm Isom}({\mathbb Q\mathbb U})$ in which every element is listed infinitely often. Let also $\tom=A_0\subseteq A_1\subseteq A_2\subseteq \ldots$ be an increasing exhaustive sequence of finite subsets of ${\mathbb Q\mathbb U}$. We shall construct an increasing sequence of finite groups $G_0\leqslant G_1\leqslant G_2\leqslant\ldots$ and finite subsets $B_0\subseteq B_1\subseteq B_2\subseteq\ldots\subseteq {\mathbb Q\mathbb U}$, where
each $G_n$ acts faithfully by isometries on $B_n$ in such a way that the action of $G_{n+1}$ on $B_{n+1}$ extends the action of $G_n$ on $B_n$, $A_n\subseteq B_n$, and for every $n$ there is some $g\in G_n$ such that $g|_{A_n}=f_n|_{A_n}$. In this case, $G=\bigcup_nG_n$ naturally acts by isometries on ${{\mathbb Q\mathbb U}}=\bigcup_n B_n$, and letting $\Gamma$ be the image of $G$ by the natural homomorphism into ${\rm Isom}({\mathbb Q\mathbb U})$, we see that $\Gamma$ is locally finite and dense in ${\rm Isom}({\mathbb Q\mathbb U})$.

To see how this is done, suppose $B_n$ and $G_n$ are defined and consider $f_{n+1}$ and $A_{n+1}$. Arbitrarily extend the action of $G_n$ on the finite subset $B_n\subseteq {\mathbb Q\mathbb U}$ to an action by isometries on ${\mathbb Q\mathbb U}$. So letting the generator of the infinite cyclic group $\Z$ act as $f_{n+1}$, we  see that this defines an action of $G_n*\Z$ by isometries on ${\mathbb Q\mathbb U}$. As $G_n$ is finite, $G_n$ has property (RZ) and so does $\Z$. By the main result of Coulbois \cite{coulbois}, also $G_n*\Z$ has property (RZ). It follows by Theorem \ref{metric} that there is a finite rational metric space $Y$ on which $G_n*\Z$ acts by isometries and an isometric  $G_n\cup\{f_{n+1}\}$-embedding of $B_n\cup A_{n+1}$ into $Y$. By the extension property of ${\mathbb Q\mathbb U}$, we can suppose that actually $B_n\subseteq B_n\cup A_{n+1}\subseteq Y\subseteq {\mathbb Q\mathbb U}$. Now let $B_{n+1}=Y$ and let $H$ be the image of $G_n*\Z$ in ${\rm Isom}(Y)$. Since $G_n$ acts faithfully on $B_n\subseteq Y$, the canonical homomorphism is injective on $G_n$ and so we can suppose that $G_n\leqslant H$. Set $G_{n+1}=H$.
\end{proof}

By exactly the same argument as above, we see that when equipped with the permutation group topology, ${\rm Aut}(\bf R)$ has a locally finite dense subgroup. This gives an alternative proof of the main result of M. Bhattacharjee and D. Macpherson \cite{macpherson}.

We see that the above method applies fairly generally to automorphism groups of many random relational structures, the main issue being whether they admit an analogue of Theorem \ref{metric}.

\end{document}